\documentclass[final,letterpaper,onecolumn]{siamltex1213}

\usepackage{algorithm}
\usepackage{algpseudocode}
\usepackage{amsmath}
\usepackage{color}
\usepackage{graphicx}
\usepackage{hyperref}
\usepackage{listings}
\usepackage{multirow}
\usepackage{textcomp}
\usepackage{verbatim}
\usepackage{subfig}

\hypersetup{pdftitle = {SpAMM in the Strong Scaling Limit}}
\hypersetup{pdfauthor = {Nicolas Bock}}
\hypersetup{pdfsubject = {Sparse Approximate Matrix Multiply (SpAMM) can be
found at https://github.com/FreeON/spammpack}}

\title{
  Solvers for $\mathcal{O} (N)$ Electronic Structure in the Strong Scaling Limit
}

\author{Nicolas Bock \and Matt Challacombe\thanks{Theoretical Division, Los
  Alamos National Laboratory, Los Alamos, NM 87544, USA,
  nicolasbock@freeon.org} \and
  Laxmikant V. Kal\'{e}\thanks{Parallel Programming Laboratory, Department of
  Computer Science, University of Illinois at Urbana-Champaign}}

\date{\today}

\begin{document}

\newcommand{\charmpp}{Charm{}\tt{++}}
\newcommand{\cpp}{C{}\tt{++}}

\maketitle

% vim: spell:spelllang=en_us:syntax=tex:nocindent:noautoindent

\begin{abstract}
  We present a hybrid OpenMP/{\charmpp} framework for solving the
  $\mathcal{O} (N)$ Self-Consistent-Field eigenvalue problem with
  parallelism in the strong scaling regime, $P\gg{N}$, where $P$ is
  the number of cores, and $N$ a measure of system size,
  \emph{i.e.}~the number of matrix rows/columns, basis functions,
  atoms, molecules, \emph{etc}. This result is achieved with a nested
  approach to Spectral Projection and the Sparse Approximate Matrix
  Multiply [Bock and Challacombe, \emph{SIAM J.~Sci.~Comput.}~{\bf 35}
    C72, 2013], and involves a recursive, task-parallel algorithm,
  often employed by generalized $N$-Body solvers, to occlusion and
  culling of negligible products in the case of matrices with decay.
  Employing classic technologies associated with generalized $N$-Body
  solvers, including over-decomposition, recursive task parallelism,
  orderings that preserve locality, and persistence-based load
  balancing, we obtain scaling beyond hundreds of cores per molecule
  for small water clusters ([H${}_2$O]${}_N$, $N \in \{ 30, 90, 150
  \}$, $P/N \approx \{ 819, 273, 164 \}$) and find support for an
  increasingly strong scalability with increasing system size $N$.
\end{abstract}



% vim: syntax=tex:spell:spelllang=en_us:nocindent:noautoindent

\section{\label{sec:01.introduction} Introduction}
% vim: syntax=tex:spell:spelllang=en_us:nocindent:noautoindent

\emph{Ab initio} electronic structure methods for the Self-Consistent-Field
(SCF) problem, involving pure density functional theory (DFT)
\cite{PhysRev.136.B864, PhysRev.140.A1133} or hybrid functionals that also
include the Fock exchange \cite{B3LYP:Becke}, offer predictive power
at low cost, finding broad utility in chemistry, biology, materials
science and drug design.  With conventional methods, solving the SCF
eigenvalue problem significantly contributes to the total
computational cost due to its steep $\mathcal{O} (N^{3})$
scaling \cite{Bylaska2009, chow2015parallel} which in practice
restricts problems to systems with $\sim{1,000}$ atoms even on large
computers \cite{quantumEspresso, Hafner2008, Valiev20101477,
gamess}. Recently, alternative methods which are $\mathcal{O} (N)$
(linear scaling) have been developed that exploit the local quantum
nature of non-metallic electronic interactions.  Early approaches to
linear scaling solutions of the SCF eigenproblem sought to exploit
this quantum locality by avoiding the pair-wise support of local basis
functions beyond a cutoff radius, leading to matrix sparsity and an
$\mathcal{O} (N)$ computational effort through iterative algorithms
based on the sparse matrix-matrix multiply (SpMM) \cite{Li1993,
Bowler2001, Bowler2002, Skylaris2005, ONETEP:2010}.  Later,
incomplete/inexact methods based on the dropping of small elements
(radial cutoffs/filtering) were developed \cite{Li1993, Daw1993,
Millam1997, Challacombe:2000:SpMM}.  While current linear scaling
methods can access systems involving $\sim{1,000,000}$
atoms \cite{Bowler2010, Miyazaki2014, VandeVondele2012}, they have yet
to enjoy widespread scientific use at scale, perhaps because the
demands of configurational sampling likewise increase with system
size. Thus, parallel algorithms that reduce the time to
solution \emph{per atom} are key in unlocking the scientific potential
of $\mathcal{O} (N)$ methods. For an excellent review and current
state of the art see Bowler \emph{et al.}~\cite{Bowler2010,
Bowler2012, Bowler2014:comment}.

A parallel SpMM implementation was first mentioned by Goringe \emph{et
al.}~as part of the CONQUEST code using a one dimensional, row-wise
matrix decomposition, although details were not
reported \cite{Goringe1997}.  Later, one of us introduced the
distributed blocked compressed sparse row (DBCSR) data format and
corresponding algorithm for distributed sparse matrix multiplication,
with space filling curve ordering and a one dimensional, row-wise
matrix decomposition based on the greedy bin packing problem,
demonstrating parallel efficiency of the SCF eigenproblem up to 128
cores
\cite{Challacombe:2000:SpMM}.  More recently, space filling curve ordering
schemes to improve locality and data layout in radial cutoff schemes
\cite{Bowler2001, Brazdova2008} and two-dimensional matrix decompositions
\cite{Borstnik201447} have lead to improved efficiencies. Bowler \emph{et
al.}~reported a scalable SpMM on 196,000 cores involving $\sim
1,000,000$ atoms \cite{Bowler2010, Miyazaki2014}, while
VandeVondele \emph{et al.}~demonstrated scalability of $\sim
1,000,000$ atoms on 46,656 cores \cite{VandeVondele2012}. Also,
generic methods for the SpMM have been developed by Bulu\c{c} \emph{et
al.}~where matrix row and columns are randomly permuted to achieve an
even load distribution, yielding high efficiencies
\cite{Buluc:2008:SpMM, Buluc2008, buluc2011parallel,
Buluc2010, 4536313}.  This approach has been adopted for quantum
chemistry with a slightly modified Cannon
algorithm \cite{cannon1969cellular}, radial cutoffs, and static load
balancing based on a fixed graph \cite{Borstnik201447}.

These parallel approaches to the $\mathcal{O} (N)$ SCF eigenvalue
problem, based on one- or two-dimensional strategies for matrix
decomposition, have established scalability in the weak regime,
$P/N \approx$ constant, where $P$ is the number of cores and $N$ is
the system size (see for example Fig.~1 of
Ref.~\cite{Bowler2014:comment}).  However, bounding communication
costs to achieve scalability beyond the weak regime remains
challenging \cite{Ballard2013} and will gain in importance for the
increasingly large, asynchronous, and heterogeneous next generation of
high performance computing systems with $P > 1,000,000$
cores\footnote{Already, the number of cores in the current top 5
supercomputers is close to or even exceeds this number: Tianhe-2 --
3,120,000 cores, Titan -- 560,640 cores, Sequoia -- 1,572,864 cores, K
-- 705,024 cores.}.  In addition, current randomization
strategies \cite{Buluc:2008:SpMM, VandeVondele2012} that forgo
locality are throttled to $\mathcal{O} (\log{P})$ \cite{Ballard2013}
due to the cost of their communication
algorithms, \emph{e.g.}~SUMMA \cite{van1997summa}; lowering these
communication costs will require either an \emph{a priori} knowledge
of sparsity patterns, or pre-computing and packing of non-zero
elements before communication \cite{Ballard2013}.  So far, even
prototypes of either strategy have yet to appear.

Recently we have developed an $N$-Body approach to the linear algebra
of data-local matrices with decay, involving the recursive occlusion
of sub-multiplicative norms based on the Cauchy-Schwarz
inequality \cite{Challacombe2010, Bock2013, spammpack}.  Besides wide
application in physical simulation \cite{WarrenGordonBell1992,
WarrenGordonBell1997, WarrenGordonBell1997b, KawaiGordonBell1999,
IshiyamaGordonBell2012, RahimianGordonBell2010, HamadaGordonBell2009,
GlosliGordonBell2007, NarumiGordonBell2006, Warren:1992:HOT,
Warren2013}, $N$-Body methods find broad applicability in statistical
learning \cite{Gray:2001:NBody, ram2010ltaps, moore2000npt,
lee2009mcmm, lee2006fgt2} and database
operations \cite{Mishra:1992:JPR:128762.128764, 4115860}.  Our Sparse
Approximate Matrix Multiply (SpAMM) algorithm is loosely comparable to
the solution of Poisson's equation through $N$-Body simulation with
radial cutoff, which has been shown recently to exhibit communication
optimal bounds, $\mathcal{O} (1/P)$, for locality preserving spatial
decompositions \cite{Driscoll2013}.  With heuristic schemes that
parlay quantum locality into spatial and temporal data locality,
together with persistence based load-balancing and three-dimensional
over-decomposition strategies, the communication cost of SpAMM may be
limited in a similar fashion.

In modern electronic structure theory there are typically four
additional ``fast'' solvers beyond the SCF eigenproblem that must
interoperate with each other, representing a tightly coupled
collective of advanced numerical methods.  Historically, these solvers
have been developed and optimized independently, involving differing
data structures and programming models (\emph{e.g.}~3-D~FFT, CSR based
SpMM, transformations of basis function to numerical
grids \emph{etc}.).  In the strong scaling regime, such a piecemeal
collection may: (a) disrupt data locality with redistributions and
transformations, (b) significantly raise the barrier to entry and
innovation, (c) exceed the ability of advanced runtime systems to load
balance multiple programming models, (d) lead to divergent rates of
error accumulation, and (e) impede deployment for trends such as fine
grained check-pointing \cite{Raicu2011, zheng2012scalable,
skarlatos2013towards},
fault-tolerance \cite{Walters:2007:FTH:1414854}, energy aware load
balancing \cite{ma2012data, tempFTSC2013} and job malleability
\cite{Kale2002, Sarood2014}.

We have recast all five solvers at the hybrid HF/DFT level of SCF
theory within the generalized $N$-Body solvers framework, including
(1) Fock exchange \cite{challacombe2014n}, (2) spectral projection
(this work), (3) inverse factorization \cite{Challacombe2015}, (4)
Coulomb summation \cite{Challacombe:1997:QCTC, Challacombe:1996:QCTCa,
Challacombe:1996:QCTCb, Challacombe:Review} and (5) the exchange
correlation problem \cite{Challacombe:2000:HiCu}.  These developments
offer a potential solution to challenges (a)-(e), through a unified
approach with a proven record of performance
\cite{WarrenGordonBell1992, WarrenGordonBell1997, WarrenGordonBell1997b,
KawaiGordonBell1999, IshiyamaGordonBell2012, RahimianGordonBell2010,
HamadaGordonBell2009, GlosliGordonBell2007, NarumiGordonBell2006,
Warren:1992:HOT, Warren2013}.  In this contribution, we develop
strategies for recursive over-decomposition and persistence-based load
balancing of the SpAMM kernel \cite{Challacombe2010, Bock2013} as
employed by spectral projection, an $\mathcal{O}(N)$ alternative to
the SCF eigenvalue problem for matrices with
decay \cite{Niklasson:2002:Pure}.  Ultimately, generic $N$-Body
frameworks and associated parallelization strategies, explored here in
part, may lead to broad horizontal support and cohesion across entire
solver collectives, enabling access to the strong scaling regime for
complex problems such as electronic structure.

It should be pointed out that the density matrix constructed through
purification schemes, such as the method of Palser and
Manolopoulos \cite{PhysRevB.58.12704} and the SP2
method \cite{Niklasson:2002:Pure} do not retain contact with the
Hamiltonian eigenspace, exponentially accumulating numerical errors
under inexact/incomplete approximation \cite{Niklasson:2003:TRS4}.  In
addition, spectral projection solvers can not be preconditioned with
the density matrix from a previous step, \emph{e.g.}~within a
molecular dynamics or structure optimization procedure, negatively
impacting overall performance \cite{PhysRevB.58.12704}.  On the other
hand, variational approaches such as the methods of Li, Nunes,
Vanderbilt \cite{Li1993}, and Daw \cite{Daw1993} retain contact with
the eigenspace of the Hamiltonian through the
gradient \cite{Bowler1999}, however, convergence can be very slow.
More recently Newton-Schulz techniques have been developed which yield
accelerated rates of convergence, and maintain direct contact with the
Hamiltonian eigenspace \cite{chen2014, Challacombe2015}.

This paper is organized as follows: In Sec.~\ref{sec:02.approach} we
describe in detail the SpAMM algorithm and in
Sec.~\ref{sec:03.task_parallel} its parallel implementations within
OpenMP and the {\charmpp} runtime.  In Sec.~\ref{sec:04.results} we
detail our methodology and show parallel scaling results for quantum
mechanical matrices with decay and demonstrate scalable
high-performance in the strong scaling limit.  Finally, we discuss our
results in Sec.~\ref{sec:05.conclusions}.

\section{\label{sec:02.approach} The Sparse Approximate Matrix Multiply}
% vim: syntax=tex:spell:spelllang=en_us:nocindent:noautoindent

A wide class of problems exist that involve matrices with decay, often
corresponding to matrix functions \cite{Benzi1999}, notably the matrix
inverse \cite{demkomosssmith, Benzi:2000:Inv}, the matrix exponential
\cite{Iserles:2000:Decay}, and in the case of electronic structure theory, the
Heaviside step function (spectral projector) \cite{McWeeny12061956,
PhysRevB.58.12704, Challacombe:1999:DMM, Challacombe:2000:SpMM,
Benzi:2007:Decay, BenziDecay2012}.  A matrix $A$ is said to decay when
its matrix elements decrease exponentially, as $\left| a_{i j} \right|
< c \,\, \lambda^{\left| i-j \right|}$, or algebraically as $\left|
a_{i j} \right| < c / ( \left| i-j \right|^{\lambda}+1 )$ with
separation $|i-j|$.  See Benzi for an excellent
discussion \cite{Benzi1999, Benzi:2000:Inv, Benzi:2007:Decay,
Benzi2013}.  In simple cases, the separation $|i-j|$ may correspond to
an underlying physical distance $|\vec{r}_{i}
- \vec{r}_{j}|$, \emph{e.g.}~of basis functions, finite
elements, \emph{etc.}~\cite{BenziDecay2012}, leading often to a strong
diagonal dominance when ordered
carefully \cite{Challacombe:2000:SpMM}.  For simple decay, truncation
in the two-dimensional vector space, \emph{e.g.}~via radial cutoff
$a_{ij} = 0$ if $| \vec{r}_{i} - \vec{r}_{j} | >
r_{\mathrm{cut}}$ \cite{Bowler1999}, a numerical threshold $a_{ij} =
0$ if $|a_{ij}| < \epsilon$ \cite{Challacombe:2000:SpMM}, or by
restricting matrix operations to a known sparsity
pattern \cite{Borstnik201447}, together with the use of a conventional
SpMM algorithm \cite{Gustavson:1978:TFA:355791.355796}, yields a
reduced complexity kernel for the iterative construction of matrix
functions.  When matrix operations are restricted to a known sparsity
pattern, the matrices retain their sparsity by construction throughout
the iterative process.  But when radial or numerical truncation
schemes are employed the matrices will fill-in unless repeatedly
filtered \cite{Bowler2012}.

\begin{figure}
  \includegraphics[width=1.0\columnwidth]{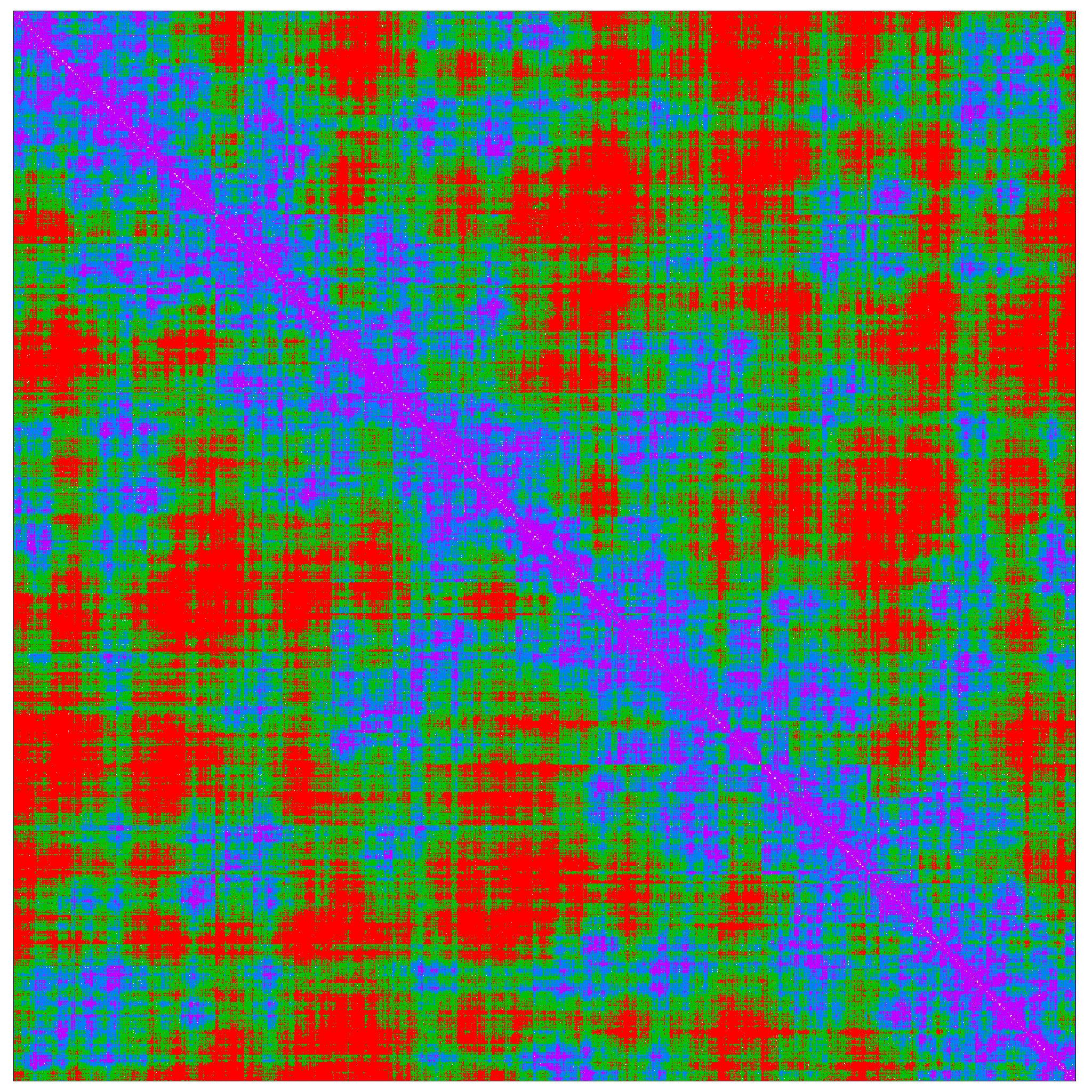}
  \caption{
    \label{fig:sparsity_pattern}
    The decay of matrix element magnitudes of a converged spectral projector
    (density matrix) for a $\left( \mathrm{H}_{2} \mathrm{O} \right)_{300}$
    water cluster at the RHF/6-31G${}^{**}$ level of theory ($n = 7500)$,
    where the molecular geometry has been reordered with a space filling
    Hilbert curve. The different colors indicate different matrix element
    magnitudes; red: $\left[ 0, 10^{-8} \right)$; green: $\left[ 10^{-8},
    10^{-6} \right)$; blue: $\left[ 10^{-6}, 10^{-2} \right)$; violet: $\left[
    10^{-2}, 1 \right]$, corresponding to approximate exponential decay.
  }
\end{figure}

Truncation may not be the most efficient or accurate approach to
exploiting decay, which can be oscillatory, involve quantum beats, or
even long range charge transfer as in the case of excited states, see
for example Ref.~\cite{Challacombe2014} and references therein.  In
addition, exploiting the secondary ``lensing'' effects in higher
dimensional operation spaces within truncation schemes is
challenging \cite{Challacombe2015}.  These effects are shown in
Fig.~\ref{fig:sparsity_pattern}, which shows a density matrix for a
large water cluster with the underlying basis ordered to preserve
locality; note the large \emph{anti-diagonal} beats, as well as the
strong clustering and segregation of elements with like magnitude.
For this type of structured matrix with non-trivial decay, the
quadtree \cite{Samet:1990:DAS:77589, Samet:2006:DBDS,
Wise:1984:RMQ:1089389.1089398, Wise:Ahnentafel}
\begin{equation}
  A^{t} = \left(
    \begin{array}{cc}
      A^{t+1}_{11} & A^{t+1}_{12} \\
      A^{t+1}_{21} & A^{t+1}_{22}
    \end{array}
  \right),
\end{equation}
where $t$ denotes the tier, pioneered in linear algebra by
Wise \emph{et.~al} \cite{Wise1990282,
springerlink:10.1007/3-540-51084-2_9, Adams:2006:SOS,
Frens:1997:AMT:263767.263789, Gottschling:2007:RMA:1274971.1274989,
Lorton:2006:ABL:1166133.1166134, Wise:2001:LSM:568014.379559},
provides a powerful framework for recursive database operations such
as the {\em metric-query} \cite{Amossen:2009:SpMMeqJoin, 4115860,
Jacox:2003:ISJ:937598.937600}, involving the lookup of sub-blocks by
magnitude, $\Vert A^t_{ij} \Vert$.  In this work we use the Frobenius
norm, which is cheap to hierarchically compute from submatrix norms,
\begin{equation}
  \Vert A \Vert = \sqrt{\sum_{ij} |A_{ij}|^{2}}
\end{equation}
and
\begin{equation}
  \Vert A^{t} \Vert = \sqrt{\sum_{i, j = 1}^{2} \Vert A^{t+1}_{ij} \Vert^{2}}.
\end{equation}
Based on this framework, the SpAMM algorithm \cite{Challacombe2010,
Bock2013}
\begin{equation}
  \label{eq:spamm_condition}
  C^t_{ij} \leftarrow C^t_{ij} + \sum_{k = 1}^{2} \left\{
    \begin{array}{lll}
      A^t_{ik} B^t_{kj}
      & \hspace{0.2cm}
      & \mbox{for} \,\, \Vert A^t_{ik} \Vert \Vert B^t_{kj} \Vert > \tau \\[0.15cm]
      0 & & \mbox{otherwise}
    \end{array}
    \right.
\end{equation}
exploits decay recursively in the three-dimensional convolution space,
with adaptive culling and occlusion of insignificant products at each
tier $t$, determined by application of the sub-multiplicative norm
inequality, $\Vert A \, B \Vert \leq \Vert A \Vert \Vert B \Vert$, and
a numerical threshold $\tau$ controlling precision.

While the discussion has so far involved dense matrices, SpAMM is
applicable to sparse matrices as well.  Also, even with dense
matrices, large values of $\tau$ correspond to an implicit truncation
and potentially a sparse product.  Relative to conventional row-column
approaches to the SpMM, the SpAMM algorithm applied to structured,
data-local matrices with decay may achieve: ({\em i}) additional
flexibility in the three-dimensional task space for domain
decomposition and load balancing (this work), ({\em ii}) the recursive
accumulation of terms with like magnitude and an $\mathcal{O} (N \lg
N) $ error accumulation \cite{Bock2013, Bini:1980:FastMM}, ({\em iii})
occlusions that occur early in recursion enabling communication
optimal approaches, ({\em iv}) a more efficient use of high level
memory chunking for message passing and low level blocking strategies
for acceleration (also this work), and ({\em v}) additional
flexibility for achieving error control within a culled volume, and
complexity reduction via lensing.

\section{\label{sec:03.task_parallel} Task Over-Decomposition}
% vim: syntax=tex:spell:spelllang=en_us:nocindent:noautoindent

One of the strengths of the generalized $N$-Body framework is that
there are many ways to realize over-decomposition on a range of
hardware, \emph{e.g.}~from long pipe GRAPE single instruction,
multiple data (SIMD) accelerators \cite{NarumiGordonBell2006,
makino1997grape, makino1998scientific, makino2003grape,
kawai1999grape} to conventional symmetric mulitprocessing (SMP) and
multiple instruction, multiple data (MIMD)
architectures \cite{4115860, Lieberman:2008:FSJ:1546682.1547260,
Warren2013, WarrenGordonBell1997, Warren:1995:HOT,
WarrenGordonBell1992, WarrenGordonBell1997b, GlosliGordonBell2007,
HamadaGordonBell2009, IshiyamaGordonBell2012, KawaiGordonBell1999,
RahimianGordonBell2010}.  Ideally, an architecture independent runtime
system seamlessly enables the recursive generation of lightweight
tasks, as OpenMP 3.0 does for SMP. However, while this feature is a
target of the Dynamic Parallelism framework of NVIDIA's CUDA
5.0 \cite{cuda} and at least partially included in a number of
parallel runtimes such as Intel's Threading Building Blocks
(TBB) \cite{IntelTBB}, Concurrent Collections (\emph{e.g.}~Intel's
CnC \cite{cnc}), Wool \cite{Wool}, Nanos++ \cite{nanox},
OpenUH \cite{openuh}, Intel's Cilk Plus \cite{IntelCilkPlus}, the Open
Community Runtime (OCR) \cite{OCR}, OpenCL \cite{opencl},
TASCEL \cite{tascel}, \emph{etc.}, full support for recursive task
parallelism is mostly unrealized for distributed memory systems at
present.  In this work, we consider simple methods for achieving
recursive task parallelism with SpAMM for the ubiquitous ``cluster of
SMP nodes'' architecture \cite{Sterling95beowulf, hargrove2001,
Chang:2011:2016802, Dayde2004} using two runtimes, OpenMP and
{\charmpp}, within a hybrid approach.

There are two main considerations in our scheme that involve memory
and task management: First, the na\"{\i}ve use of task parallelism at
the SMP level, with either OpenMP or {\charmpp}, has the potential to
involve non-contiguous memory and high packing/unpacking overheads
when redistributing memory between nodes, potentially negatively
impacting overall performance.  Additionally, the cache hierarchy of
modern CPUs with small, local caches and large, shared last level
caches should not be ignored.  Second, an explicitly allocated,
unrolled octree is a necessary structure that enables {\charmpp} to
manage tasks involving occlusion and culling as well as node-level SMP
work due to limitations of the load-balancing framework implemented in
{\charmpp}.

Thus, we allocate contiguous chunks of size $N_{c} \times N_{c}$ to
hold a full sub-quadtree together with a $N_{b} \times N_{b}$ blocking
at the lowest level.  The chunks are processed using OpenMP and the
code can potentially be used without modification on the Intel Xeon
Phi coprocessor and through automatic source code
translation \cite{noaje2011source, Martinez2011, sabne2012effects,
Chang2013} on GPGPUs.  In addition, the use of OpenMP removes
{\charmpp} compile and runtime dependencies for single-node
applications, potentially significantly simplifying the build process.

We expect the overall performance to be influenced by several
competing size-dependent effects: (1) The ratios $N/N_{c}$ and
$N_{c}/N_{b}$ limit the maximum number of tasks available for
load-balancing for {\charmpp} and OpenMP respectively, and (2) the
leaf node size $N_{b}$ affects the performance of memory access
through the CPU's cache hierarchy and the potential for vectorization
and convolution space compression.  While we previously demonstrated
that a highly specialized and optimized dense kernel can lead to
competitive performance for very small dense submatrices of $N_{b} =
4$ \cite{Bock2013}, the use of manually tuned assembly code renders
this approach rigid with respect to submatrix granularity and width of
SIMD vectors.  Thus, in this work, we implemented a simple kernel with
three nested loops and leave low-level optimizations to the compiler.

\subsection{\label{sec:OpenMP} OpenMP}

Shown in Alg.~\ref{alg:spamm_OpenMP} is the SMP parallel
implementation within the OpenMP application programming
interface; \texttt{SpAMM\_omp} recursively walks a transient octree
generated dynamically on the stack through the OpenMP 3.0 tasking
feature \cite{openmp}.  Guided by the binary convolution of matrix
quadtrees $A^t$ and $B^t$ at each tier $t$ (line 3), the implicit
octree traversal may be sparse and irregular due to culling and
occlusion (line 4) based on the sub-multiplicative matrix norm
inequality, $\Vert A \, B \Vert \leq \Vert A \Vert \Vert B \Vert$.
The parallel tree traversal is extended through untied
OpenMP \texttt{task}s (line 5) and recursive calls
to \texttt{SpAMM\_omp} (line 6).  Per node synchronization (to ensure
appropriate variable lifetimes) is achieved through the OpenMP
\texttt{taskwait} statement (line 9).  Finally, at the leaf tier, dense matrix
products are performed and the result reduced into the $C$ quadtree
(line 12), with a data race on $C$ prevented through explicit use of
OpenMP locks (lines 11, 13). While other approaches to address data
write contention are certainly possible, \emph{e.g.}~OpenMP reductions
or atomics, we found good on-node parallel scaling using explicit
locks, and defer such potentially performance enhancing details to a
forthcoming article.

In this work we have made only modest effort to optimize the {\tt
SpAMM\_omp} implementation, or even the dense contraction on line 12.
In Ref.~\cite{Bock2013} we showed that accuracies better than the
native GEMM are possible also with $N$-scaling, but that difficult,
platform specific optimizations were necessary; we are currently
developing a corresponding OpenMP algorithm and are investigating the
use of compiler vectorization and OpenMP 4.0 SIMD
constructs \cite{openmp}.

\begin{algorithm}
  \caption{
    \label{alg:spamm_OpenMP}
    The OpenMP SpAMM algorithm, recursively multiplying matrices $C
    \gets A \times B$ under a SpAMM tolerance $\tau$.  The function
    matrix arguments are pointers to tree nodes.  }
  \begin{algorithmic}[1]
    \Function{SpAMM\_omp}{$\tau$, $t$, $A^t$, $B^t$, $C^t$}
    \If{${t}<{}$ depth}
    \ForAll{$\left\{ i, j, k \,\, \middle| \,\, C^t_{ij} \gets A^t_{ik} B^t_{kj} \right\}$}
    \If{$\left\| A_{ik} \right\| \left\| B_{kj} \right\| > \tau$}
    \Comment{Culling}
    \State{OpenMP {\tt task untied}}
    \State{\Call{SpAMM\_omp}{$\tau$, t+1, $A^{t+1}_{ik}$, $B^{t+1}_{kj}$, $C^{t+1}_{ij}$}}
    \EndIf
    \EndFor
    \State{OpenMP {\tt taskwait}}
    \Else
    \State{\Call{omp\_set\_lock}{}}
    \Comment{Acquire OpenMP lock on $C$}
    \State{$C \gets C + A \times B$}
    \Comment{Dense product}
    \State{\Call{omp\_unset\_lock}{}}
    \Comment{Release OpenMP lock on $C$}
    \EndIf
    \EndFunction
  \end{algorithmic}
\end{algorithm}

\subsection{{\charmpp}}\label{charmpp}

{\charmpp} \cite{charm} is a mature runtime environment on distributed
memory platforms available for all major supercomputer systems,
allowing for efficient scalable high performance
implementations \cite{Menon2013, tempFTSC2013, CharmAppsBook:2013,
LifflanderPLDI, 2007_ChaNGaScaling, cosmo2007, NamdSC11,
NAMDEncyEntry10, NamdSC07, Zhao2013}.  In the message-driven execution
model of {\charmpp}, code and data are encapsulated in {\cpp} objects
called ``chares'' which are initially placed by static load balancing
algorithms.  Dynamic persistence-based load balancing strategies
migrate chares transparently during solver execution based on load and
communication measurements from previous solver iterations and
efficiently optimize load distribution and communication cost. The
{\charmpp} runtime transparently manages chare placement and migration
and proxy objects are used to send messages to particular chare
instances or groups thereof without explicit specification of their
location.  Chares can be grouped in multi-dimensional sparse arrays or
used as ``singleton'' objects.

Persistence-based load balancing exploits temporal and spatial
localities in iterative solvers through decomposition of the load and
communication graph.  Since the dynamic load balancing strategies of
{\charmpp} only consider chares organized in arrays persistently
instantiated across solver operations and load balancing, the
multiplication octree has to be explicitly stored in memory (as
opposed to the transient stack based ``storage'' used in the SMP
implementation).  Note that such persistent allocation of the
multiplication octree could aid efficient load balancing across
molecular dynamics or structure optimization steps,
see \emph{e.g.}~the impressive scaling of astrophysics
applications \cite{cosmo2007, jetley2008massively, Warren2013}.  The
nodes of the matrix quadtree between root, $t = 0$, and chunks, $t =
t_{c}$, given by $N_{c}$, are stored in a stack of two-dimensional
chare arrays of size $2^{t} \times 2^{t}$ each.  The corresponding
unrolled octree is stored in three-dimensional chare arrays with
occlusion and culling carried out iteratively, tier-by-tier, until the
chunk level at which {\tt SpAMM\_omp} is invoked.

Data and work locality are exploited through the communication aware
load balancing strategies in {\charmpp}.  However, at the time of this
writing, a bug in the {\charmpp} runtime \cite{charmbug_445} prevents
the use of \emph{sparse} load balanced chare arrays. As a work-around,
we mark chares that correspond to pruned tree nodes with a boolean
data member, {\tt isDisabled == true}, introducing a $\mathcal{O}
(N^{3})$ communication component with a prefactor found to be
negligible.

\begin{algorithm}
  \caption{
    \label{alg:spamm_charm}
    The SpAMM algorithm in the {\charmpp} programming language. Tree
    occlusion is done by iterating over the three-dimensional
    multiplication chare arrays, \textsc{convolution[$d$]}. In
    {\charmpp} a call such as \textsc{convolution[$t$].occlude}
    translates into a broadcast to all array elements
    of \textsc{convolution[$t$]}.
  }
  \begin{algorithmic}[1]
    \Function{SpAMM\_charm}{$\tau$, $A$, $B$, $C$}
    \For{$t \ge 0 \land t < d$}
    \State{\Call{convolution[$t$].occlude}{$\tau$}}
    \Comment{See Alg.~\ref{alg:spamm_charm_culling}}
    \EndFor
    \State{\Call{convolution[$d$].multiply}{}}
    \State{\Call{convolution[$d$].store}{}}
    \EndFunction
  \end{algorithmic}
\end{algorithm}

The {\charmpp} algorithm is outlined in Alg.~\ref{alg:spamm_charm} and
proceeds in three phases. In the first phase, the multiplication
octree is constructed iteratively over the top tiers of the
three-dimensional chare arrays, shown in lines 2 and 3 of
Alg.~\ref{alg:spamm_charm}.  This phase is a breadth-first
implementation of the SpAMM algorithm and retains the full complexity
reduction of the depth-first, recursive implementation,
Alg.~\ref{alg:spamm_OpenMP}. In each iteration of this phase, a
broadcast message is sent to all multiplication chares of tier $t$
(line 3) executing the \textsc{occlude} method on the enabled array
elements, shown in Alg.~\ref{alg:spamm_charm_culling}. The scalar
products of the eight matrix norms of the $A$ and $B$ nodes of the
next tier are formed, lines 6-10 of
Alg.~\ref{alg:spamm_charm_culling}, and Eq.~\ref{eq:spamm_condition}
is used to decided whether to enable or disable the corresponding
multiplication chares.  Disabled multiply chares ({\tt isDisabled ==
true}) are skipped during the next iteration of the pruning phase,
shown in lines 2-4 of Alg.~\ref{alg:spamm_charm}.

\begin{algorithm}
  \caption{
    \label{alg:spamm_charm_culling}
    Tree occlusion in the {\charmpp} programming language of the
    multiplication chare element on tier $t$ with index $(i, j,
    k)$. In {\charmpp} the call \textsc{convolution[$t+1$]($i, j,
    k$).enable} translates into a direct message to
    the \textsc{enable} method of the multiplication chare element
    with index ($i, j, k$) on tier $t+1$.
  }
  \begin{algorithmic}[1]
    \Function{occlude}{$\tau$}
    \Comment{On tier $t$, index $(i, j, k) \in [1, 2^{t}]$}
    \If{isDisabled}
    \State{return}
    \EndIf
    \ForAll{$\left\{ i', j', k' \,\, \middle| \,\, C^{t+1}_{i', j'} \gets A^{t+1}_{i' k'} B^{t+1}_{k' j'} \right\}$}
    \If{$\Vert A^{t+1}_{i' k'} \Vert \Vert B^{t+1}_{k' j'} \Vert > \tau$}
    \State{\Call{convolution[$t + 1$]($i', j', k'$).enable}{}}
    \Else
    \State{\Call{convolution[$t + 1$]($i', j', k'$).disable}{}}
    \EndIf
    \EndFor
    \EndFunction
  \end{algorithmic}
\end{algorithm}

During the second phase, line 5 of Alg.~\ref{alg:spamm_charm}, the SMP
SpAMM code is called to compute the $N_{c} \times N_{c}$ submatrix
products in each remaining, enabled multiplication chare, and the
results are stored in a temporary variable local to the chare.  In the
final phase, line 6 of Alg.~\ref{alg:spamm_charm}, all temporary
matrix products are gathered in the
\textsc{store} method, summed, and added to the corresponding chares of $C$.
Since the {\charmpp} runtime guarantees exclusive execution of chare
instances, explicit locking or other means of synchronization as in
the OpenMP implementation are not necessary.

\subsection{The OpenMP/Charm++ Hybrid}

In our hybrid approach, we found the best performance with one
{\charmpp} Processing-Element (PE) per node, OpenMP commanding all
on-node threads and $N_{c} \times N_{c}$ quadtree chunking as
discussed above.  This approach avoids the problem of packing and
unpacking fragmented memory during chare migration, enabling use of a
single \texttt{memcpy}, which is efficient in standard libraries such
as \texttt{libc}.  Certainly, optimal chunk and block sizes are likely
to be application dependent, an issue beyond the scope of the current
work.  A further complication of the hybrid approach involves the
issue of local \emph{vs.}~absolute addressing; by wrapping an address
offset with convenience macros, the OpenMP application programming
interface given in Alg.~\ref{alg:spamm_OpenMP} can be used without
modification.

\section{\label{sec:04.results} Results}
% vim: spell:spelllang=en_us:syntax=tex:nocindent:noautoindent

In this work we consider scalability of the SpAMM kernel in the
context of spectral projection \cite{McWeeny12061956,
PhysRevB.58.12704, Challacombe:1999:DMM, Challacombe:2000:SpMM,
Benzi:2007:Decay, BenziDecay2012}, an alternative to explicitly
solving the SCF eigenvalue problem \cite{szabo1996modern}.  Spectral
projection involves nested construction of the matrix Heaviside
step-function from the effective SCF Hamiltonian (Fockian), in our
case computed in a basis of atom-centered
functions \cite{szabo1996modern}.  In this work, tightly converged,
dense matrices for a sequence of water clusters were computed at the
B3LYP/6-31G** level of theory \cite{B3LYP:Becke}
using \texttt{FreeON}, a suite of programs for $\mathcal{O} (N)$
quantum chemistry \cite{FreeON}.  This sequence of water clusters
corresponds to standard temperature and pressure, and has been used in
a number of previous studies \cite{Challacombe:Review,
Challacombe:1996:QCTCb, Challacombe:1997:QCTC, burant:8969,
ESchwegler97, millam:5569, daniels:425, ochsenfeld:1663, Bock2013}.
The 6-31G** basis set introduces 5 basis function per hydrogen atom
and 15 basis functions per oxygen atom, yielding 25 basis functions
(and matrix rows and columns) per water molecule.  A key aspect of
this work is ordering of the atom indices with the locality preserving
Hilbert curve, Ref.~\cite{Challacombe:2000:SpMM} and references
therein, yielding clustering and segregation of elements by magnitude
as in Fig.~\ref{fig:sparsity_pattern}.

In a previous study \cite{Bock2013}, we reported linear scaling
computational complexities and SpAMM errors as the max norm of the
difference between the SpAMM product and a dense reference product.
Here, we consider SpAMM errors that accumulate in iterative
application of the second order spectral projection scheme
(SP2) \cite{Niklasson:2002:Pure}.  In all cases, the SP2 solver was
run to convergence, taking 40 iterations.  Values of $\tau = 10^{-6},
10^{-8}$, and $10^{-10}$ are considered for scaling experiments, with
$\tau = 10^{-6}$ corresponding to extreme truncation (a highly sparse
representation).

All OpenMP scaling studies were run on a fully allocated 48-core,
4-socket AMD Opteron 6168 (Magny Cours architecture) system running at
1.9 GHz using \texttt{GNU gcc} 4.6.3, and a 24 core, 2-socket AMD
Opteron 6176 (Magny Cours architecture) system running at 2.3 GHz
using \texttt{GNU gcc} 4.7.2.  The {\charmpp} scaling studies were run
on the largest open computer cluster at Los Alamos National Laboratory
(LANL), ``mustang'', which consists of 1,600 dual socket AMD Opteron
6176 (Magny Cours) nodes for a total of 38,400 cores using \texttt{GNU
gcc} 4.7.2.  All tests used the \texttt{-O2} level of compiler
optimization.

\subsection{Error Accumulation}

\begin{figure}
  \includegraphics[width=1.0\columnwidth]{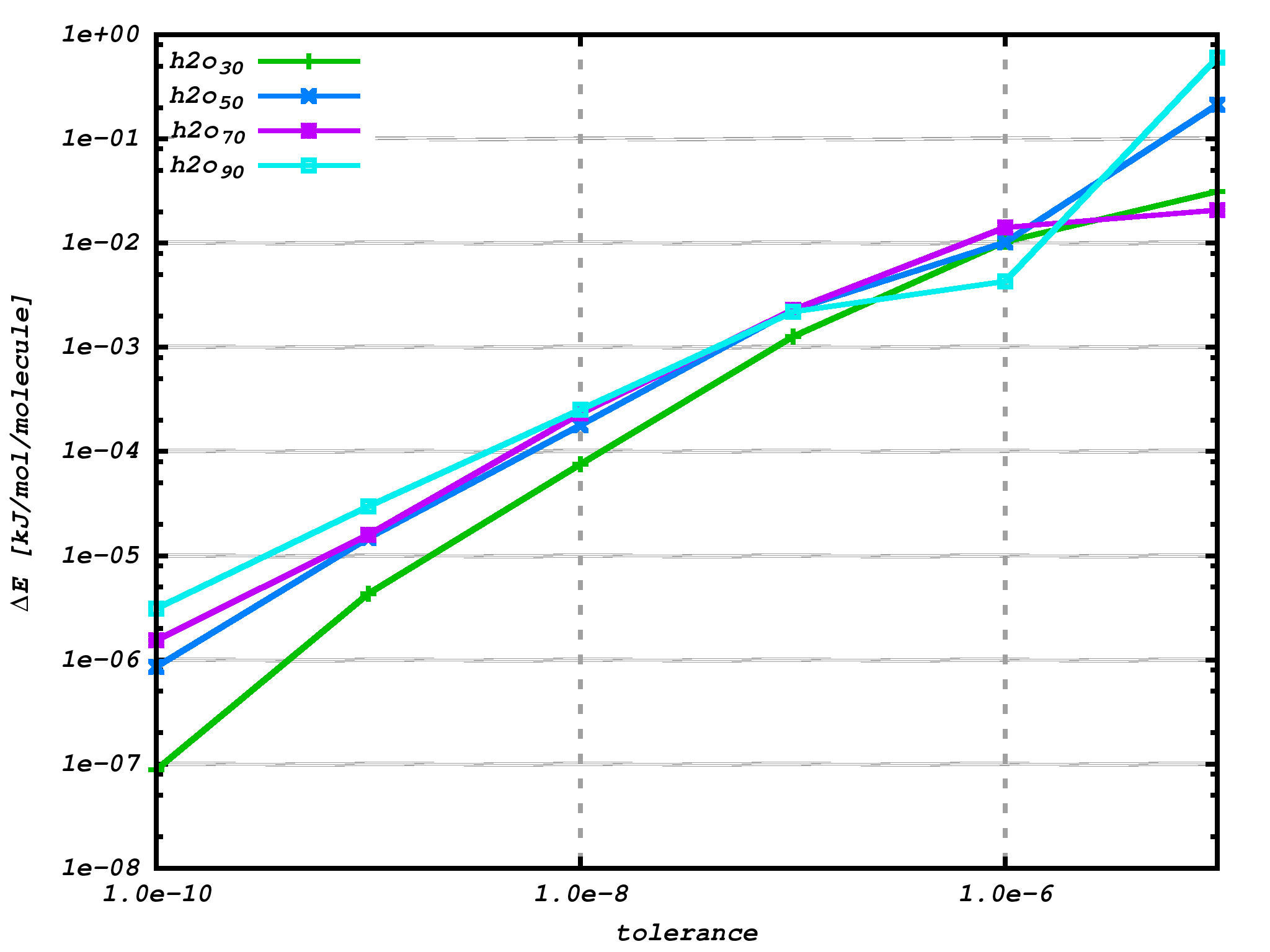}
  \caption{
    \label{fig:absolute_error}
    The absolute error of the energy after 40 iterations of the
    spectral projection method for different water clusters in
    B3LYP/6-31G${}^{**}$.
    }
\end{figure}

The accumulation of error in spectral projection due to the SpAMM
kernel is computed here as $\mathrm{Tr} [ F \, (P - \tilde{P}) ]$,
where $F$ is the Fockian, $\tilde{P}$ is the approximate density
matrix computed from $F$ with $\tau \ne 0$, and $P$ is a reference
computed with $\tau=0$.  These errors are reported in
Fig.~\ref{fig:absolute_error}, demonstrating that the error per
molecule exhibits no significant system size dependence for the cases
studied here, in agreement with our earlier results on the max norm
error behavior of SpAMM, Figs.~5.2 and 5.3 of Ref.~\cite{Bock2013}.
Roughly, these results suggest that chemical accuracy (1 kcal/mol or
4.184 kJ/mol \cite{Heil1970}) may be retained with $10^{5}$ water
molecules and a SpAMM threshold of $\tau = 10^{-10}$.  As discussed in
Sec.~\ref{sec:02.approach}, the control of accumulated errors in the
spectral projection solver is challenging due to the non-variational
nature of the solver.  However, our results indicate good error
control even under extreme truncation conditions ($\tau = 10^{-6}$)
due to the recursive occlusion and culling based on the
sub-multiplicative norm inequality, as opposed to matrix element
truncation directly in the vector space.

\subsection{OpenMP scaling}

\begin{figure}
\centering{
\subfloat[\texttt{SpAMM\_omp} restricted to 1 thread.]{
\includegraphics[width=0.5\columnwidth]{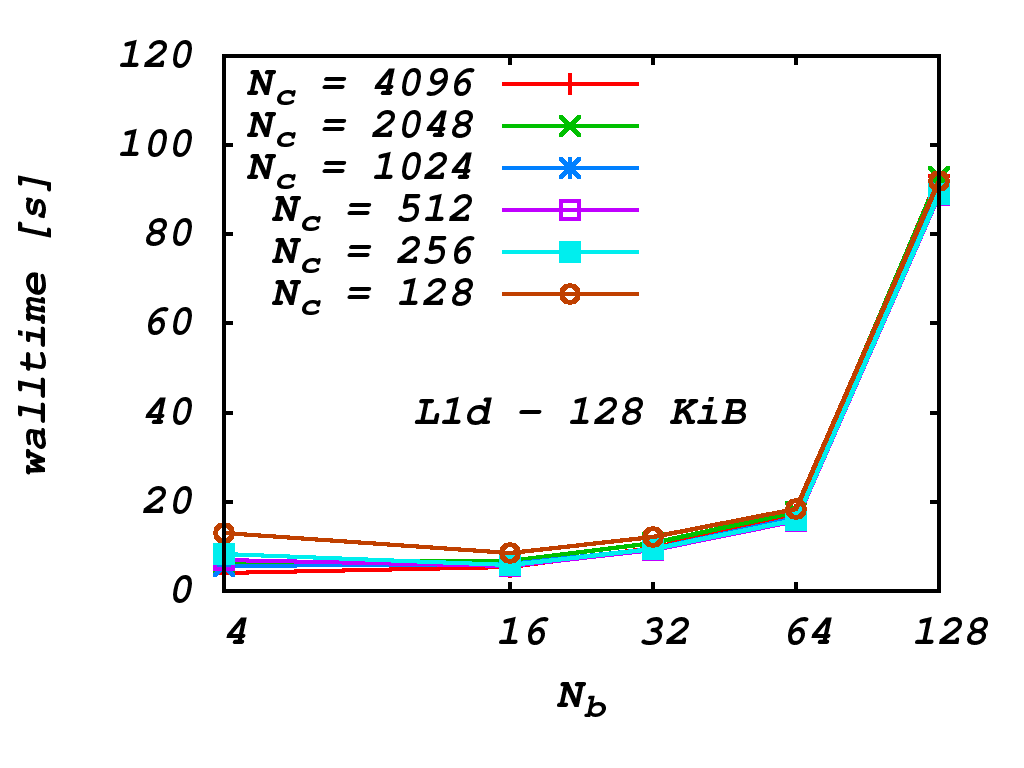}
\hspace{-0.4cm}
}
\subfloat[\texttt{SpAMM\_omp} on all 48 threads.]{
\includegraphics[width=0.5\columnwidth]{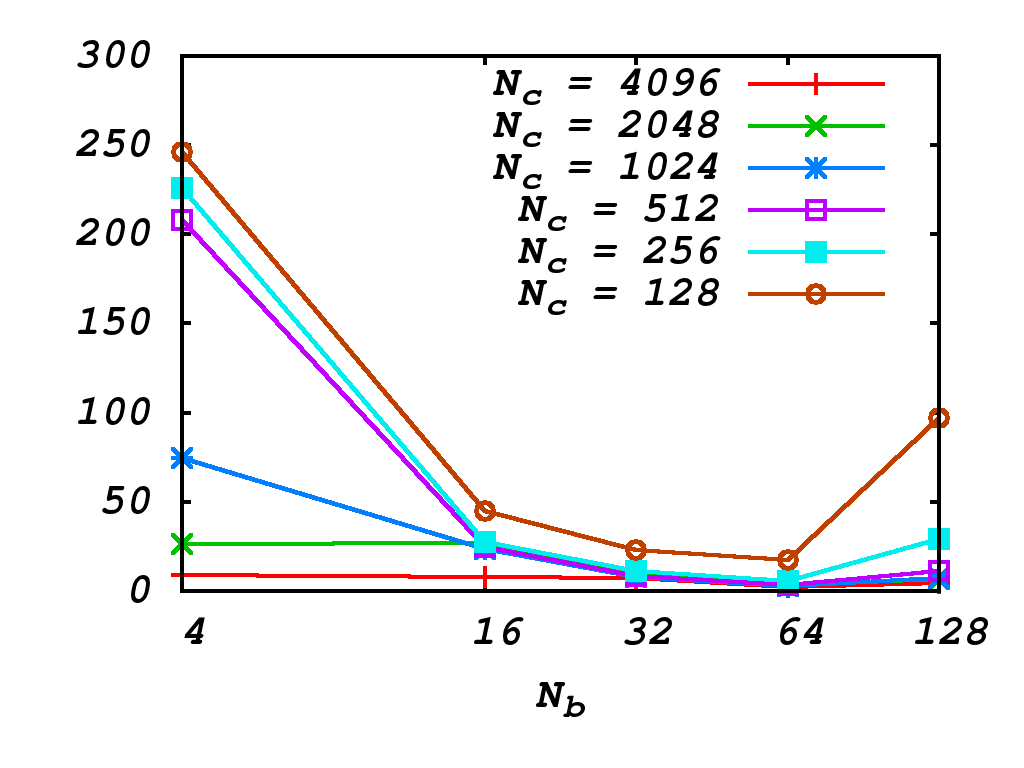}
}
}
\caption{\label{fig:openmp_scaling}
Total time of matrix product for (H${}_{2}$O)${}_{90}$
of \texttt{SpAMM\_omp} with $\tau = 10^{-6}$ restricted to one thread,
(a), and on 48 threads, (b), on the Opteron 6168.  In serial, we find
the measured walltime to be independent of $N_{c}$, but to depend
strongly on $N_{b}$.  Note that $N_{b} > 64$ exceeds L1d leading to
significant performance loss.  An increase in compression due to
smaller granularity leads to decreasing walltime with decreasing
$N_{b}$, and an optimal block size of $N_{b} = 16$.  Tests with $\tau
= 0$ indicate that the optimal block size in serial without
compression is $N_{b} = 64$.  On 48 threads, the shortest walltime
shifts from $N_{b} = 16$ to $N_{b} = 64$, shown in (b), which
indicates poor memory access performance and is potentially due to a
lack of cache/core affinity.  It is worth noting that we find a large
spread in performance for $N_{b} = 4$ across the chunk sizes tested.}
\end{figure}

\begin{figure}
\centering
\subfloat[48-core Opteron 6168.]{
\includegraphics[width=0.5\columnwidth]{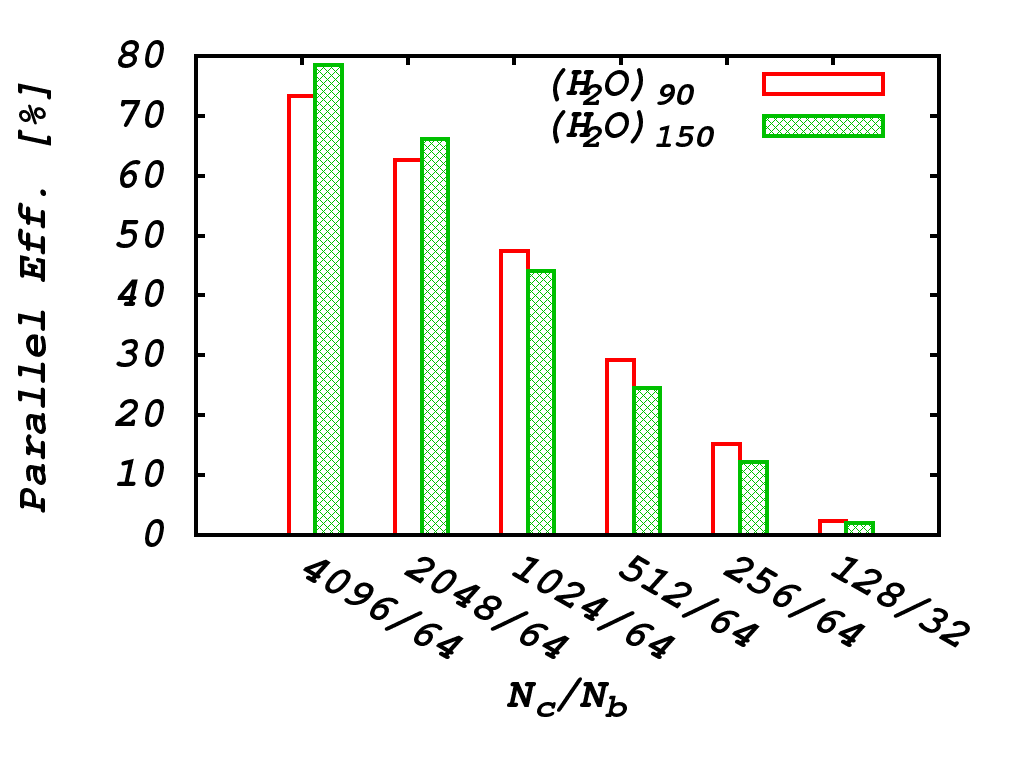}
\hspace{-0cm}
}
\subfloat[24-core Opteron 6176.]{
\includegraphics[width=0.5\columnwidth]{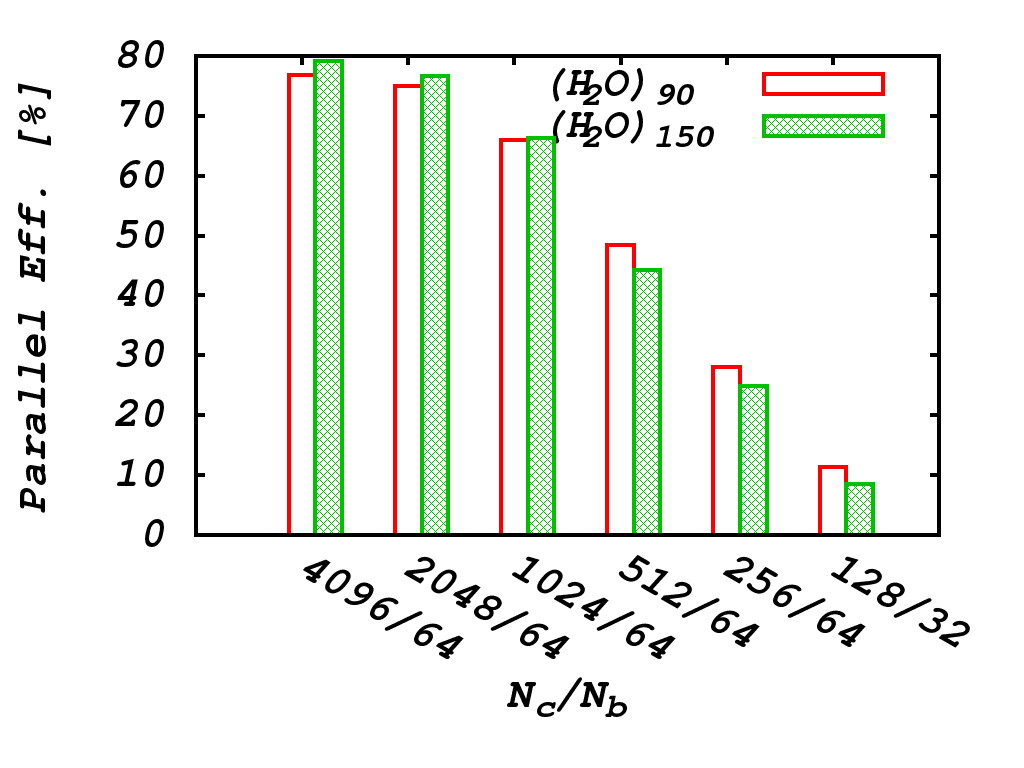}
}
\caption{\label{fig:openmp_parallel_efficiency}
Parallel efficiency of OpenMP code with different granularities for
B3LYP/6-31G** (H${}_{2}$O)${}_{90}$ and (H${}_{2}$O)${}_{150}$ with
$\tau = 10^{-6}$. The best performing combination of $N_{c}/N_{b}$ at
the maximum number of threads ($P = 48$) was chosen for each value of
$N_{c}$ tested. Note that the native matrix size of $2,250 \times
2,250$ is padded to $4,096 \times 4,096$.  There is little difference
between the two water clusters.  While larger chunks exhibit good
parallel efficiency, we find a significant drop of parallel scaling
for the smaller chunks.}
\end{figure}

\texttt{SpAMM\_omp} (Alg.~\ref{alg:spamm_OpenMP}) was benchmarked for the last
SP2 iteration of the (H${}_{2}$O)${}_{90}$ and (H${}_{2}$O)${}_{150}$
density matrices with a loose SpAMM tolerance of $\tau = 10^{-6}$.
These examples are well within the linear scaling regime for SpAMM
calculations \cite{Bock2013}, yet small enough to probe a lower
molecule/core ratio available on modern SMP platforms.  In addition,
such a loose SpAMM tolerance leads to irregular work loads and data
access, potentially challenging the OpenMP runtime.  Shown in
Fig.~\ref{fig:openmp_scaling} are walltimes for the smaller water
cluster, (H${}_{2}$O)${}_{90}$, under several combinations of $N_{b}$
and $N_{c}$ in serial (OpenMP restricted to one thread), (a), and on
48 threads, (b).  We find the overall performance and parallel scaling
to be significantly influenced by $N_{c}$ and $N_{b}$.  While the
performance in serial, Fig.~\ref{fig:openmp_scaling} (a), is mostly
independent of the tested chunk sizes, $N_{c} = \{ 128, 256, 512,
1024, 2048, 4096 \}$, it is significantly impacted by the size of
$N_{b}$.  For large values, $N_{b} > 64$, we note a steep rise in
walltime which we attribute to a lack of temporal locality in cache
due to the size of per core L1d of 128KiB.  The improving performance
with decreasing block size, $N_{b} \leq 64$, is due to increasing
compression in convolution space due to finer granularity.  Tests with
$\tau = 0$ indicate that the optimal block size in serial is $N_{b} =
64$.  However, on 48 cores we find that the shortest walltime shifts
from $N_{b} = 16$ to $N_{b} = 64$, which indicates poor memory access
performance and might be due to a lack of thread/data affinity.  Also,
on 48 threads we find the walltime to depend more strongly on $N_{c}$,
which we attribute to a lack of potential tasks for shallow trees,
indicated by small ratios $N_{c}/N_{b}$.

In addition to walltime, it is instructive to investigate the parallel
efficiency of \texttt{SpAMM\_omp}.  For each value of $N_{c}$ we chose
the value of $N_{b}$, yielding the fastest (lowest walltime)
performance at 48 threads, and calculated the parallel efficiency as
$E(P) = T(1)/(P \,\, T(P))$, where $P$ denotes the number of threads,
shown in Fig.~\ref{fig:openmp_parallel_efficiency}.  We note that
efficiencies up to 80\% can be achieved for large chunks.  As
$N_{c}/N_{b}$ and the number of potential tasks decreases
($4096/64 \rightarrow 262,144$ potential tasks, $2048/64 \rightarrow
32,768$ potential tasks, \dots) load balancing becomes increasingly
challenging with a significant decline in parallel efficiency.
Compared to (H${}_{2}$O)${}_{90}$, the larger water cluster,
(H${}_{2}$O)${}_{150}$, exhibits a slightly superior parallel scaling
for the larger chunks.  However, the qualitative behavior with
decreasing chunk size remains the same.

\subsection{{\charmpp} Scaling}

\begin{figure}
\centering
\subfloat[First SP2 iteration, $\tau = 10^{-10}$.]{
  \includegraphics[width=0.45\columnwidth]{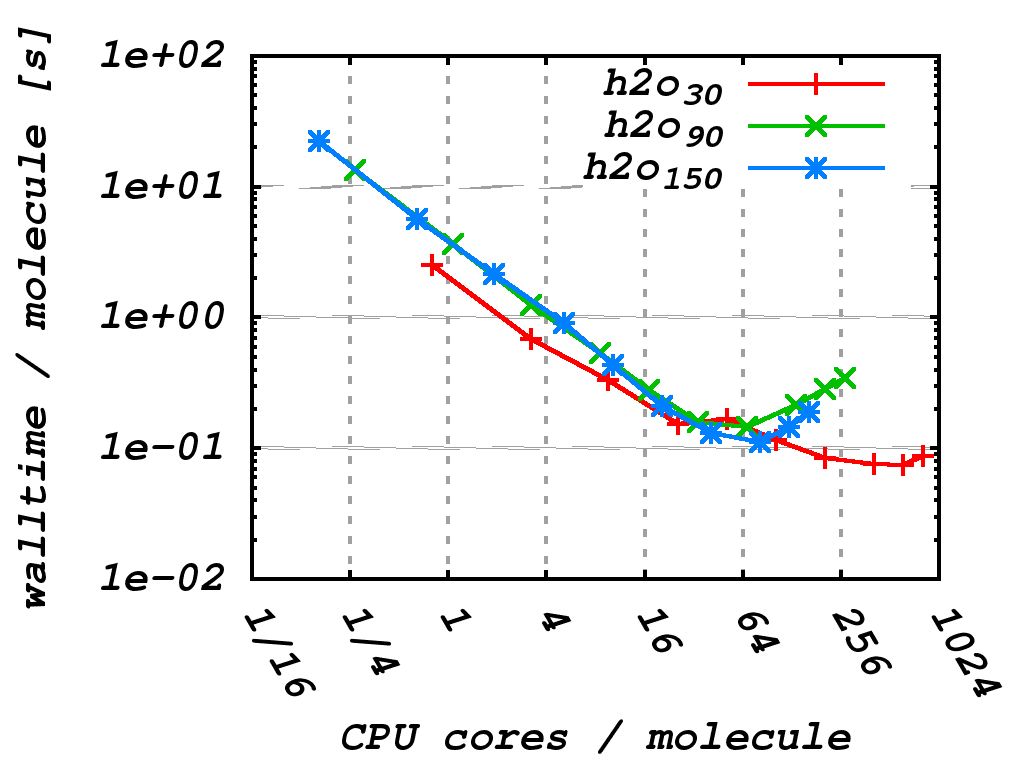}
}
\subfloat[First SP2 iteration, $\tau = 10^{-8}$.]{
  \includegraphics[width=0.45\columnwidth]{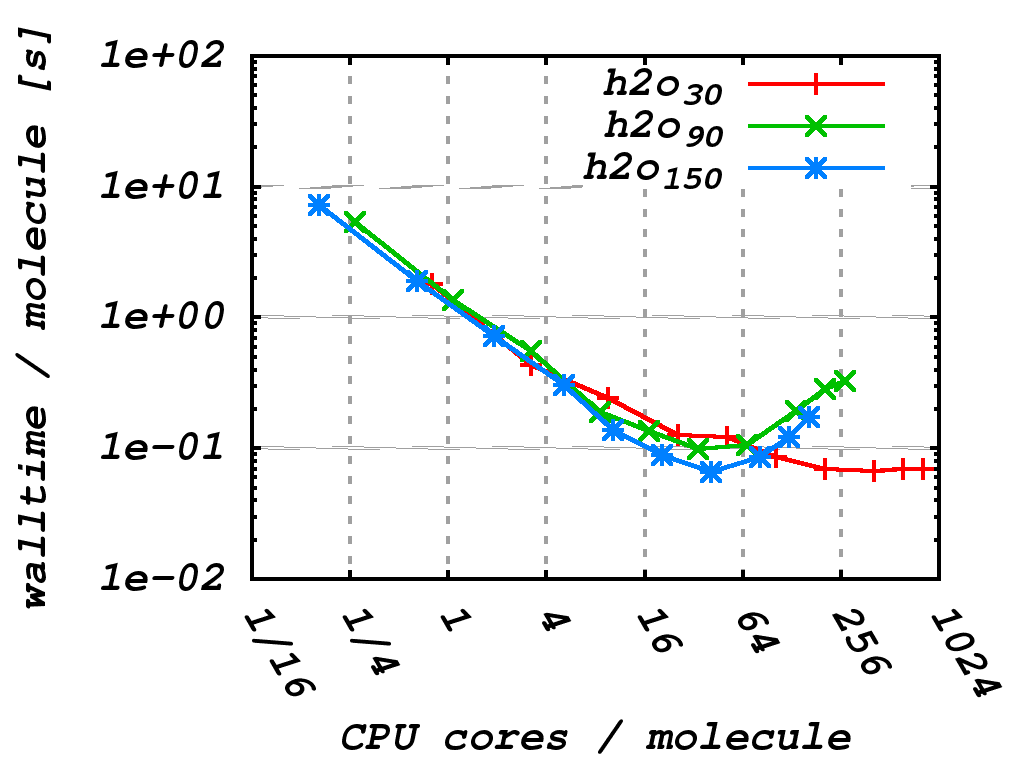}
}
\hspace{0mm}
\subfloat[First SP2 iteration, $\tau = 10^{-6}$.]{
  \includegraphics[width=0.45\columnwidth]{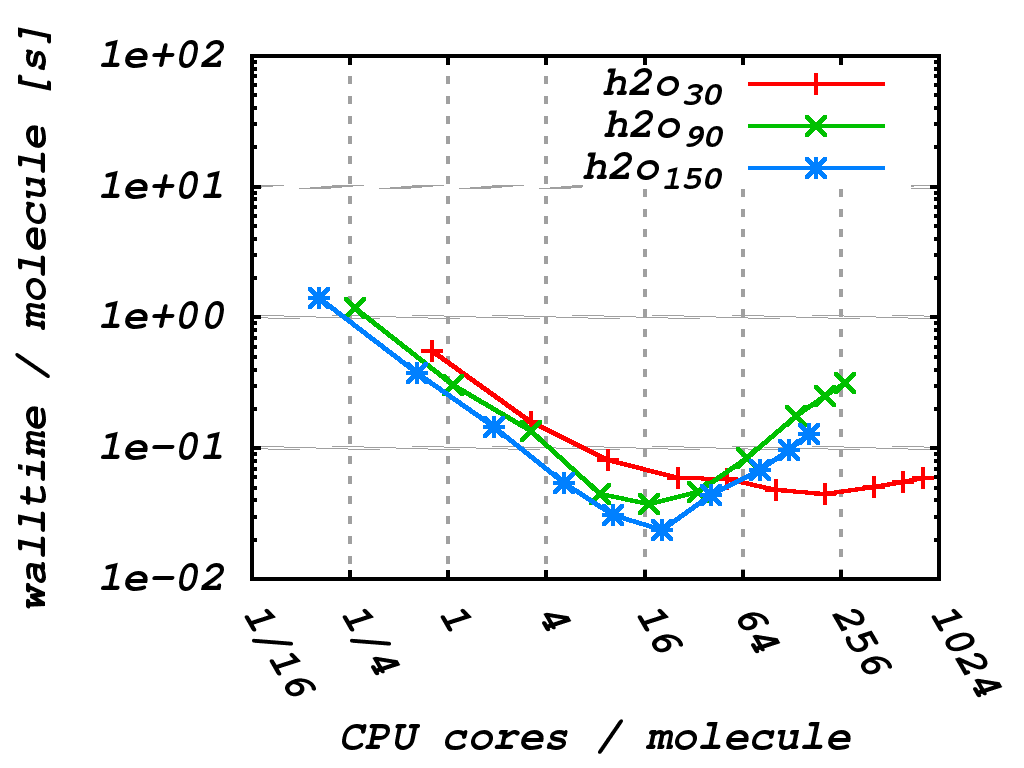}
}
\subfloat[Last SP2 iteration. $\tau = 10^{-10}$.]{
  \includegraphics[width=0.45\columnwidth]{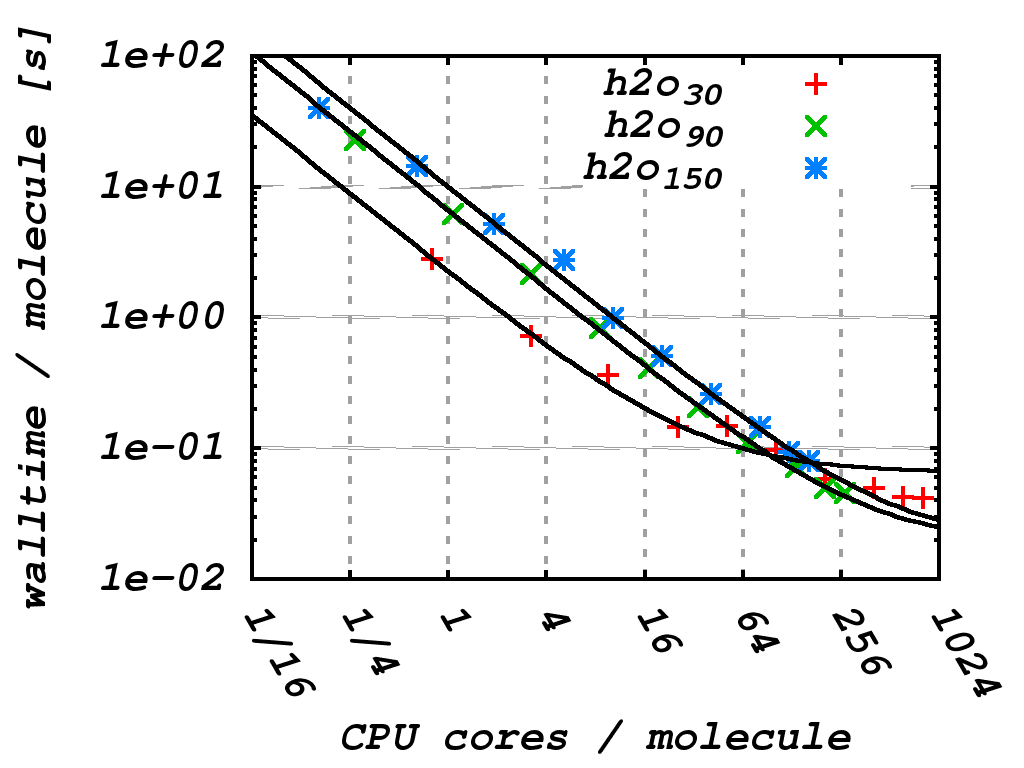}
}
\hspace{0mm}
\subfloat[Last SP2 iteration, $\tau = 10^{-8}$.]{
  \includegraphics[width=0.45\columnwidth]{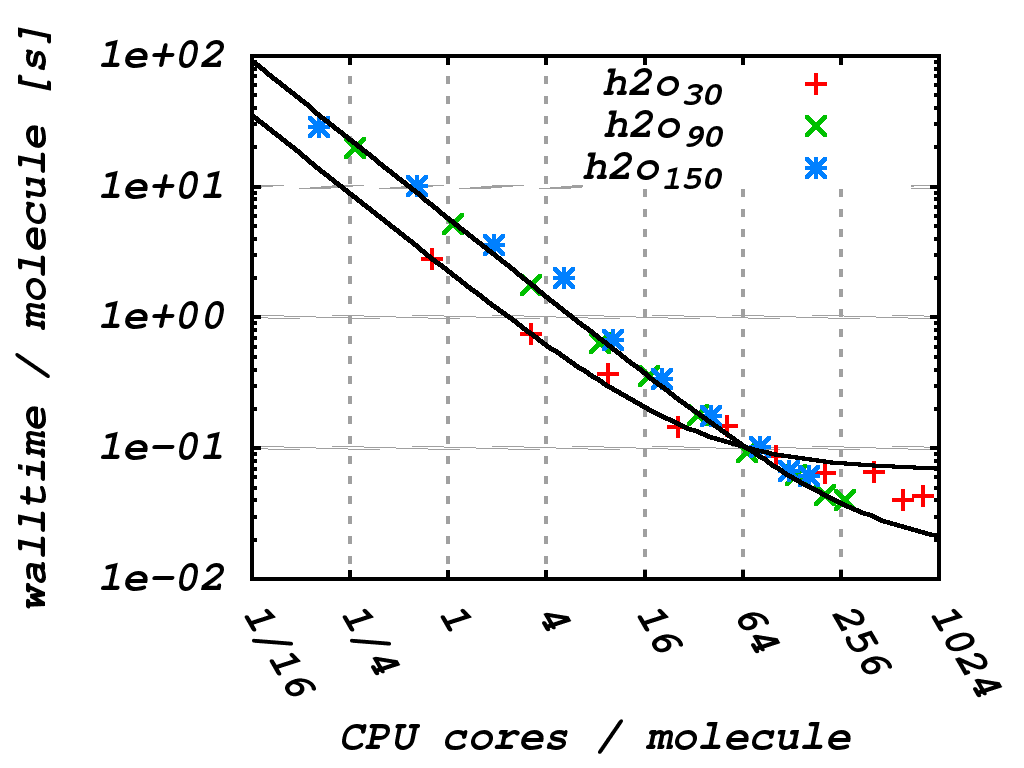}
}
\subfloat[Last SP2 iteration. $\tau = 10^{-6}$.]{
  \includegraphics[width=0.45\columnwidth]{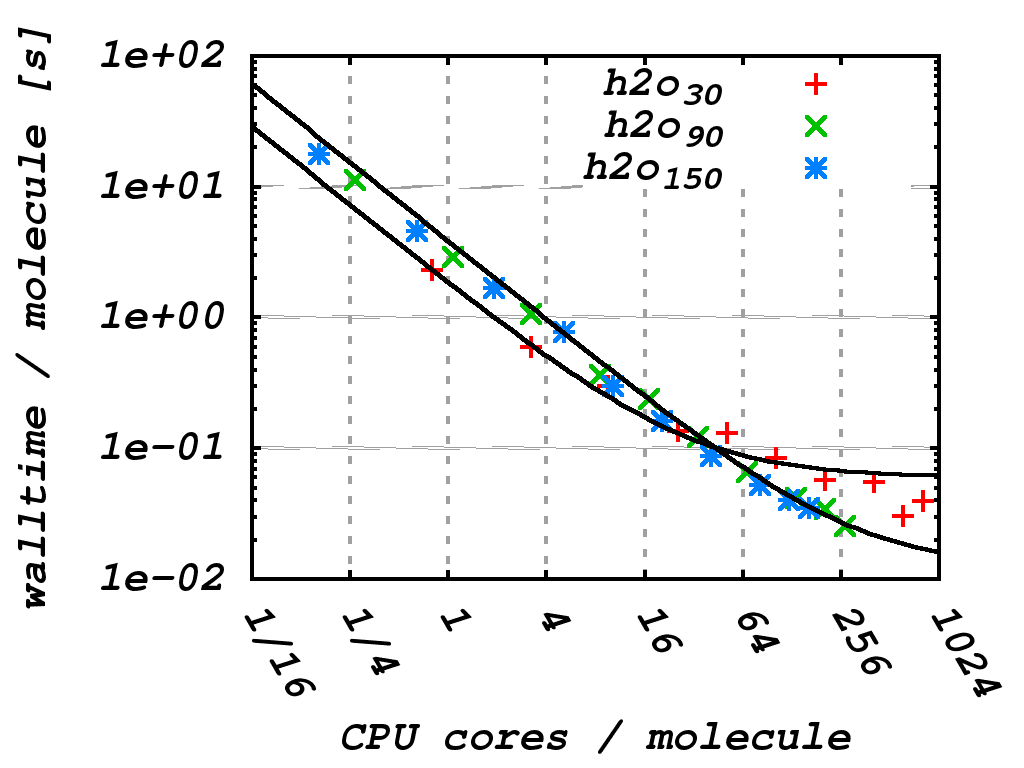}
}
\caption{
  \label{fig:charm_scaling_per_molecule}
  Shown in panels (a)-(c), scaling results of the first SP2 iteration
  under different thresholds. Shown in panels (d)-(f), scaling results
  of the last and fully load-balanced iteration (iteration 40) of SP2,
  under the same thresholds. As a guide, fits to Amdahl's law,
  $T_{s}+T_{p}/p$, are shown as solid lines, see
  Table \ref{table:charm_scaling} for fitting parameters.
}
\end{figure}

This study involved scaling with the progression $P = 24 \times 2^m$,
up to 24,576 cores (1024 nodes) on LANL's largest open computer
cluster ``mustang''.  The study consisted of spectral projection via
the SP2 method until convergence (40 iterations) with $m = 1,2, \dots,
10$, and $\tau = 10^{-6}, 10^{-8}$, and $10^{-10}$.  The initial data
distribution during the first SP2 iteration was given by the
{\charmpp} default static load balancer.  After each iteration of the
SP2 algorithm, the \texttt{GreedyCommLB} load balancer of {\charmpp}
was called to migrate matrix and multiply chares in order to rebalance
work and data. To demonstrate the efficiency of the
\texttt{GreedyCommLB} load balancer, we show walltime vs.~cores for the first
iteration (only statically balanced), panels (a)-(c) of
Fig.~\ref{fig:charm_scaling_per_molecule}, and for the final
iteration, panels (d)-(f) of
Fig.~\ref{fig:charm_scaling_per_molecule}.  Notice that the difference
in wall time between the first and last iteration is due to matrix
fill-in (the decay slows from Fockian to density matrix).  On the
first iteration, we observe scaling roughly to $P = 30 \,\, N$,
corresponding to the default {\charmpp} data distribution.  After a
few iterations however, the communication aware
persistence-based \texttt{GreedyCommLB} load balancer dynamically
migrates chares to achieve a balance of very high quality. In
applications, the persistence-based load balance will remain effective
between SCF cycles, and also as atomic-positions gradually evolve,
{\em e.g.}~in a molecular dynamics simulation, geometry optimization,
{\em etc.}, mitigating inefficiencies associated with the first
iterations.

In Table~\ref{table:charm_scaling}, we list parameters for the fits to
Amdahl's law, $T^{\,\tau}_{s} + T^{\,\tau}_{p}/p$, corresponding to
the fitted lines in panels (d)-(f) of
Fig.~\ref{fig:charm_scaling_per_molecule}.  Also given in
Table~\ref{table:charm_scaling} are the corresponding break-even core
counts, $P^{\,\tau}_{\mathrm{even}} = T^{\,\tau}_{p}/T^{\,\tau}_{s}$,
the ratio between parallel and serial components.  The break-even core
count is a conservative estimate of the core count at which additional
scaling becomes ineffective due the left-over serial component, which
was found to be 1-3 seconds in all cases.  Also, we notice a
pronounced decrease in the parallel component with increasing values
of $\tau$, due to sparse-irregular effects.  It should be noted that
this analysis is a useful quantitative guide despite it being
simplistic in ignoring more subtle scaling effects such as for example
the scaling behavior of communication collectives and the details of
network topology.

\begin{table}
  \begin{center}
    \begin{tabular}{l|r|r|r|r|r}
      \hline
      \hline
      $\tau$
      & (H${}_{2}$O)${}_{N}$
      & matrix
      & $T_{s}$ [s]
      & $T_{p}$ [s]
      & $P_{\mathrm{even}}$ \\
      \hline
      \multirow{3}{*}{$10^{-10}$}
      &  30 &   750 & 1.95 &   1,981 &  1,016 \\
      &  90 & 2,250 & 1.66 &  53,481 & 32,252 \\
      & 150 & 3,750 & 2.74 & 224,825 & 81,909 \\
      \hline
      \multirow{3}{*}{$10^{-8}$}
      &  30 &   750 & 2.05 &   1,981 &   966  \\
      &  90 & 2,250 & 1.41 &  46,459 & 32,898 \\
      & 150 & 3,750 & 2.48 & 149,819 & 60,364 \\
      \hline
      \multirow{3}{*}{$10^{-6}$}
      &  30 &   750 & 1.80 &   1,621 &    899 \\
      &  90 & 2,250 & 1.11 &  30,983 & 28,007 \\
      & 150 & 3,750 & 2.55 &  65,284 & 25,553 \\
      \hline
      \hline
    \end{tabular}
  \end{center}
  \caption{
    \label{table:charm_scaling}
    Fit parameters for Amdahl's law, $T_{s} + T_{p}/P$, corresponding
    to the curves in panel (d)-(f) of
    Fig.~\ref{fig:charm_scaling_per_molecule}.  Matrix dimensions are
    shown in the third column labeled ``matrix''.  Also listed is the
    break-even core-count $P_{\mathrm{even}} = T_{p}/T_{s}$, providing
    a conservative estimate of parallel scalability.
    }
\end{table}

\section{\label{sec:05.conclusions} Conclusions}
% vim: spell:spelllang=en_us:syntax=tex:nocindent:noautoindent

Relative to the $\approx 4$ heavy atoms/core granularity achieved in
the weak limit by advanced parallel methods \cite{Bowler2014:comment},
the default ``static'' distribution of work exhibited by our
OpenMP/{\charmpp} implementation achieves roughly $P = 30 \,\, N$, as
shown in panels (a)-(c) of Fig.~\ref{fig:charm_scaling_per_molecule}.
Assuming 1 water molecule $\approx$ 2 heavy atoms, our default is
$\sim 60 \times$ more scalable.  Once persistence is employed however,
our results extend into the strong scaling regime, yielding $P =
400 \,\, N$ to $600 \,\, N$ as inferred from
Table \ref{table:charm_scaling}.  For working accuracies and larger
systems, \emph{e.g.}~$N \gg 150$ and $\tau \in \{ 10^{-8},
10^{-12} \}$, we expect substantially better results as suggested by
Table \ref{table:charm_scaling}.  We also expect substantially better
results for problems with slower decay, as for example problems
involving semi-conductors and metal oxides.

Based on the results given in Fig, the (very modest) serial component
seems to be due to the {\charmpp} runtime.  Larger calculations on
larger computers will allow the reliable collection of diagnostics, as
well as examination of the relationships between data locality and
communication.

While our recursive, depth-first implementation of SpAMM with OpenMP
exhibits good parallel scaling for larger chunk sizes, further
improvements, including parallel performance at fine granularities,
may require more explicit approaches to exploiting the inherent
temporal and spatial localities present in SpAMM and to make contact
with the deep memory cache hierarchy of the Magny-Cours architecture.
Also, other computer platforms may not exhibit the pronounced
non-uniform memory access (NUMA) effects common to the AMD Magny-Cours
architecture, and we may expect the parallel performance
of \texttt{SpAMM\_omp} on those platforms to show improved scaling.
Finally, it is known that the runtime has significant impact on the
performance of SpAMM-like workloads \cite{Podobas2015}, and other
programming frameworks might lead to improved parallel scaling.

These satisfactory results follow from over-decomposition of the
three-dimensional convolution space, relative to conventional methods
that involve decomposition in one or two dimensions, and from runtime
systems that support the irregular task parallelism inherent in the
generalized $N$-Body solvers framework.  The ability to recursively
generate singleton chares would greatly simplify the implementation of
$N$-Body methods, and enhance their efficiency by eliminating the
explicit management of tree-traversal as explained in
Section \ref{charmpp}.  This prospect, together with a unified code
base for $N$-Body solver collectives (based on established prototypes
\cite{challacombe2014n, Challacombe2015, Challacombe:1997:QCTC, Challacombe:1996:QCTCa,
Challacombe:1996:QCTCb, Challacombe:Review, Challacombe:2000:HiCu}),
may offer a simple and well posed approach to meeting the challenges
of increasing hardware complexity.

The software written for and used in this study is available online
at \url{http://www.freeon.org/spammpack} \cite{spammpack}, licensed
under the terms of the BSD license \cite{bsd}.

\section{\label{sec:06.acknowledgements} Acknowledgments}
% vim: syntax=tex:spell:spelllang=en_us:nocindent:noautoindent

NB and MC thank the LDRD program for funding this research under LDRD-ER grant
20110230ER. Both would also like to acknowledge generous support from the
Ten-Bar Caf\'{e} providing stimulating and helpful discussions. This article
was released under LA-UR-14-22050. The Los Alamos National Laboratory is
operated by Los Alamos National Security, LLC for the NNSA of the USDoE under
Contract No. DE-AC52- 06NA25396.

\bibliographystyle{siam}
\bibliography{spamm_article}

\end{document}